\documentstyle{article}
\input{amssym}
\begin{document}
\title{Lie group analysis for short pulse equation}
\author{M.~Nadjafikhah\thanks{School of Mathematics, Iran University of Science and Technology, Narmak,
Tehran 1684613114, Iran. Tel \& Fax : +98-021-73225426. e-mail: m\_nadjafikhah@iust.ac.ir}}
\maketitle
\begin{abstract}
In this paper, the classical Lie symmetry analysis and the
generalized form of Lie symmetry method are performed for a
general short pulse equation. The point, contact and local
symmetries for this equation are given. In this paper, we
generalize the results of H.~Liu and J.~Li \cite{LiuLi}, and add
some further facts, such as optimal system of Lie symmetry
subalgebras and two local symmetries.\\
{\bf Keywords:} Short pulse equation, Lie symmetry analysis, Point, contact and local symmetries.
\end{abstract}

\maketitle
\section{Introduction}
Nonlinear PDEs arising in many applied fields like the biology,
fluid mechanics, plasma physics and optics, systems of impulse and
neural networks, etc, and exhibit a rich variety of nonlinear
phenomena. The investigation of the exact solutions plays an
important role in the study of nonlinear systems. In this paper,
we find Lie point symmetries, third order local symmetries,
optimal system of these two type symmetries, and corresponding
invariant solutions for a general short pulse equation:
\begin{eqnarray}\label{SPE}
\textrm{SPE}\;:\;u_{xt}=\alpha u+\frac{1}{3}\beta (u^3)_{xx}
\end{eqnarray}
where $u=u(x,t)$ is the unknown real function and subscripts
denote differentiation w.r.t. $x$ and $t$; $\alpha$ and $\beta$
are nonzero real parameters.

This general SPE was derived by T. Schafer and C.E. Wayne
\cite[p.94]{ScWa} as a model equation describing the propagation
of ultra-short light pulses in silica optical fibres. In
\cite{ScWa,SaSa}, many results are obtained about the special SPE:
\begin{eqnarray}\label{SPE-1}
u_{xt}=u+\frac{1}{6}(u^3)_{xx}.
\end{eqnarray}
\section{Lie contact and point symmetries}
Let $J^1=J^1({\Bbb R}^2,{\Bbb R})$ be the jet space with
coordinates $(x,t,u,u_x,u_t)$. Let
\begin{eqnarray}\label{eq:1}
{\bf v}=\xi\,\partial_x+\tau\,\partial_t+\eta\,\partial_u
+\eta^x\,\partial_{u_x}+\eta^t\,\partial_{u_t}
\end{eqnarray}
be an infinitesimal Lie contact symmetry of (\ref{SPE}), where
$\xi$, $\tau$ and $\eta$ are functions $J^1({\Bbb R}^2,{\Bbb
R})\to{\Bbb R}$, and
\begin{eqnarray}\label{eq:2}
Q=\xi u_x+\tau u_t-\eta
\end{eqnarray}
be the characteristic function of (\ref{eq:2}). Thus
\begin{eqnarray}\label{eq:3}
&\xi=Q_{u_x},\;\; \tau=Q_{u_t},\;\;\eta=u_xQ_{u_x}+u_tQ_{u_t}-Q,&\\
&\eta^{x}=-Q_x-u_xQ_u,\;\;\eta^{t}=-Q_t-u_tQ_u.&\nonumber
\end{eqnarray}
Then, the (\ref{eq:1}) is an infinitesimal Lie contact symmetry of
the (\ref{SPE}) if and only if
\begin{eqnarray}\label{eq:4}
{\bf v}^{(2)}\Big(u_{xt}-\alpha u-\frac{1}{3}\beta
(u^3)_{xx}\Big)=0,\qquad u_{xt}=\alpha u+\frac{1}{3}\beta
(u^3)_{xx},
\end{eqnarray}
where ${\bf v}^{(2)}$ is the second prolongation of $\bf v$:
\begin{eqnarray}\label{eq:6}
{\bf v}^{(2)}=Q\,\partial_u+{\rm D}_xQ\,\partial_{u_x}+{\rm
D}_tQ\,\partial_{u_t}+{\rm D}_x^2Q\,\partial_{u_{xx}}+{\rm
D}_x{\rm D}_tQ\,\partial_{u_{xt}}+{\rm D}_t^2Q\,\partial_{u_{tt}},
\end{eqnarray}
where ${\rm D}_x$ and ${\rm D}_t$ are total derivative w.r.t $x$
and $t$, respectively.

By subsituting $u_{x,t}$ from second equation of (\ref{eq:6}) in
the first equation, we find a polynomial of $u_{xx}$ and $u_{tt}$
with functional coefficients of $(x,t,u,u_x,u_t)$. Its
coefficients must be zero:
\begin{eqnarray}\label{eq:7}
&&\hspace{-10mm} Q_{u_x,u_x}=0,\quad Q_{u_x,u_t}=0,\quad Q_{u_t,u_t}=0,\nonumber\\
&&\hspace{-10mm} u_tQ_{u,u_x}+\alpha uQ_{u_x,u_x}+Q_{t,u_x}=0,\quad \alpha uQ_{u_t,u_t}+u_xQ_{u,u_t}+Q_{x,u_t}=0,\nonumber\\
&&\hspace{-10mm} u_x^2Q_{u,u}+2u_xQ_{x,u}+Q_{xx} +\alpha
u\big(\alpha Q_{u_t,u_t}+2u_xQ_{u,u_t}+2Q_{x,u_t}\big)=0,\\
&&\hspace{-10mm} u_tQ_{u,u_t}-5u_xQ_{u,u_x}-5Q_{x,u_x}+Q_{t,u_t}-4u Q_{u_x,u_t}-2u_xQ_u=0,\nonumber\\
&&\hspace{-10mm} u_xu_tQ_{u,u}+u_xQ_{t,u}+u_tQ_{x,u}+Q_{x,t}\nonumber\\
&&\hspace{-5mm} +\alpha u\big(u_tQ_{u,u_t}+Q_{t,u_t}+
Q_{x,u_x}+Q_u\big)+\alpha\big(u_xQ_{u_x}+u_tQ_{u_t}-Q\big)=0.\nonumber
\end{eqnarray}
After solving the determining system (\ref{eq:7}), one finds that
\begin{eqnarray}\label{eq:8}
Q=C_1\,u_x+C_2\,u_t+C_3\,(xu_x-tu_t-3u);
\end{eqnarray}
where, $C_1$, $C_2$ and $C_3$ are arbitrary constants. Therefore,
\paragraph{Theorem}
{\it The SPE (\ref{SPE}) has a $3-$dimensional Lie algebra $\goth
g$ of point symmetries, generated by the infinitesimal generators
\begin{eqnarray}\label{eq:9}
{\bf v}_1=\partial_x,\;\;{\bf v}_2=\partial_t,\;\; {\bf
v}_3=x\,\partial_x-t\,\partial_t+3u\,\partial_u,
\end{eqnarray}
and commutating table
\begin{eqnarray}\label{eq:10}
\begin{array}{ccccc}
\hline
  [\,,\,] & {\bf v}_1 & {\bf v}_2
  & {\bf v}_3 \\ \hline
  {\bf v}_1  & 0             & 0            & {\bf v}_1 \\
  {\bf v}_2  & 0             & 0            & -{\bf v}_2\\
  {\bf v}_3  & -{\bf v}_1 & {\bf v}_2 & 0\\
  \hline\end{array}\end{eqnarray}
The SPE (\ref{SPE}) has not any non-point contact symmetry.}
\section{Invariant solutions and its classification}
The one-parameter groups $G_i$ generated by the base of $\goth g$
are as follows:
\begin{eqnarray}\label{eq:11}
&&\hspace{-10mm} G_1\;:\;\exp(\varepsilon{\bf v}_1)\cdot(x,t,u)=(x+\varepsilon,t,u),\nonumber\\
&&\hspace{-10mm} G_2\;:\;\exp(\varepsilon{\bf v}_2)\cdot(x,t,u)=(x,t+\varepsilon,u),\\
&&\hspace{-10mm} G_3\;:\;\exp(\varepsilon{\bf
v}_3)\cdot(x,t,u)=\big(e^{\varepsilon}x,e^{-\varepsilon}t,e^{3\varepsilon}u\big),\nonumber
\end{eqnarray}
where $\varepsilon$ is a real number.

Since each group $G_i$ is a symmetry group of SPE (\ref{SPE}) and
if $u=f(x,y)$ is a solution of the SPE (\ref{SPE}), so are the
following functions
\begin{eqnarray}\label{eq:12}
u=f(x+\varepsilon,t),\quad u=f(x,t+\varepsilon),\quad
u=f\big(e^{\varepsilon}x,e^{-\varepsilon}t,e^{-3\varepsilon}u\big),
\end{eqnarray}
where $\varepsilon$ is an arbitrary real number. Thus, for the
arbitrary combination ${\bf v}=c_1{\bf v}_1+c_2{\bf v}_2+c_3{\bf
v}_3\in{\goth g}$, the SPE (\ref{SPE}) has the following solution:
\begin{eqnarray}\label{eq:13}
u=f\big(e^{\varepsilon_3}x+\varepsilon_1,e^{-\varepsilon_3}t+\varepsilon_2,e^{-3\varepsilon_3}u\big),
\end{eqnarray}
where $\varepsilon_i$ are arbitrary real numbers.

Let $G$ be the symmetry Lie group of SPE (\ref{SPE}). Now $G$
operates on the set of solutions $S$ of SPE (\ref{SPE}), and
$s\cdot G$ be the orbit of $s$, and $H$ be a subgroup of $G$.
Invariant $H-$solutions $s\in S$ are characterized by equality
$s\cdot S=\{s\}$. If $h\in G$ is a transformation and $s\in S$,
then
\begin{eqnarray}\label{eq:14}
h\cdot(s\cdot H)=(h\cdot s)\cdot (hHh^{-1}).
\end{eqnarray}
Consequently, every invariant $H-$solution $s$ transforms into an
invariant \break $hHh^{-1}-$solution (Proposition 3.6 of
\cite{Ovsi}). Therefore, different invariant solutions are found
from similar subgroups of $G$. Thus, classification of invariant
$H-$ solutions is reduced to the problem of classification of
subgroups of $G$, up to similarity. An optimal system of
$s-$dimensional subgroups of $G$ is a list of conjugacy
inequivalent $s-$dimensional subgroups of $G$ with the property
that any other subgroup is conjugate to precisely one subgroup in
the list. Similarly, a list of $s-$dimensional subalgebras forms
an optimal system if every $s-$dimensional subalgebra of $\goth g$
is equivalent to a unique member of the list under some element of
the adjoint representation: $\tilde{\goth h}={\rm
Ad}(g)\cdot{\goth h}$. Let $H$ and $\tilde{H}$ be connected,
$s-$dimensional Lie subgroups of the Lie group $G$ with
corresponding Lie subalgebras ${\goth h}$ and $\tilde{\goth h}$ of
the Lie algebra ${\goth g}$. Then $\tilde{H}=gHg^{-1}$ are
conjugate subgroups if and only $\tilde{\goth h}={\rm
Ad}(g)\cdot{\goth h}$ are conjugate subalgebras (Proposition 3.7
of \cite{Ovsi}). Thus, the problem of finding an optimal system of
subgroups is equivalent to that of finding an optimal system of
subalgebras, and so we concentrate on it.

For the one-dimensional subalgebras, the classification problem is
essentially the same as the problem of classifying the orbits of
the adjoint representation, since each one-dimensional subalgebra
is determined by a nonzero vector in Lie algebra symmetries of SPE
(\ref{SPE}) and so to "simplify" it as much as possible. The
adjoint action is given by the Lie series
\begin{eqnarray}\label{eq:15}
\hspace{-7mm}\mbox{Ad}(\exp(\varepsilon{\bf v}_i){\bf
v}_j)\!=\!{\bf v}_j\!-\!\varepsilon[{\bf v}_i,{\bf
v}_j]\!+\!\frac{\varepsilon^2}{2}[{\bf v}_i,[{\bf v}_i,{\bf
v}_j]]\!-\!\cdots,
\end{eqnarray}
where $i,j=1,\cdots,3$. Let $F^{\varepsilon}_i:{\goth
g}\rightarrow{\goth g}$ defined by ${\bf
v}\mapsto\mbox{Ad}(\exp(\varepsilon{\bf v}_i){\bf v})$, for
$i=1,\cdots,3$. Therefore, if ${\bf v}=c_1{\bf v}_1+c_2{\bf
v}_2+c_3{\bf v}_3\in{\goth g}$, then
\begin{eqnarray}\label{eq:16}
&&\hspace{-10mm}F^{\varepsilon_1}_i({\bf
v})=(c_1+\varepsilon_1c_3){\bf v}_1+c_2{\bf
v}_2+c_3{\bf v}_3,\nonumber\\
&&\hspace{-10mm}F^{\varepsilon_2}_i({\bf v})=c_1{\bf
v}_1+(c_2+\varepsilon_2c_3){\bf
v}_2+c_3{\bf v}_3,\\
&&\hspace{-10mm}F^{\varepsilon_3}_i({\bf
v})=e^{-\varepsilon_3}c_1{\bf v}_1+e^{\varepsilon_3}c_2{\bf
v}_2+c_3{\bf v}_3.\nonumber
\end{eqnarray}
Applying these transformations, one can show that
\paragraph{Theorem}
{\it An one-dimensional optimal system of ${\goth g}$ is
\begin{eqnarray}\label{eq:17}
{\bf v}_1+a{\bf v}_2,\quad b{\bf v}_1+{\bf v}_2,\quad {\bf v}_3,
\end{eqnarray}
where $a$ and $b$ are real constants; and, a two-dimensional
optimal system of ${\goth g}$ is given by
\begin{eqnarray}\label{eq:18}
{\bf v}_1,{\bf v}_2,\quad {\bf v}_1,{\bf v}_3,\quad {\bf v}_2,{\bf
v}_3.
\end{eqnarray}}
\section{Local symmetries of SPE}
One can generalize one-parameter Lie groups of point
transformations with infinitesimal generators in the
characteristic form ${\bf v}=Q(x,t,u,u_x,u_t)\,\partial_u$ to
one-parameter $s$-order local transformations with infinitesimal
generators of the form
\begin{eqnarray}\label{eq:19}
{\bf v}=Q(x,t,u,\partial u,\partial^2u,\cdots,\partial^s
u)\,\partial_u,
\end{eqnarray}
where the infinitesimal components depend on derivatives of $u$ up
to some finite order $s\geq1$. The prolongation of $\bf v$ is
given by
\begin{eqnarray}\label{eq:20}
&&\hspace{-14mm}{\bf
v}^{(\infty)}=Q\,\partial_u+{\rm D}_xQ\,\partial_{u_x}+{\rm D}_tQ\,\partial_{u_t}+{\rm D}_x^2Q\,\partial_{u_{xx}}\\
&&\hspace{-3mm}+{\rm D}_x{\rm D}_tQ\,\partial_{u_{xt}}+{\rm
D}_t^2Q\,\partial_{u_{tt}}+\cdots.\nonumber
\end{eqnarray}
where ${\rm D}_x$ and ${\rm D}_t$ are total derivative w.r.t $x$
and $t$, respectively \cite{BlAnCh}.

Then, for $s=3$, (\ref{eq:20}) is an infinitesimal local symmetry
of the (\ref{SPE}) if and only if
\begin{eqnarray}\label{eq:21}
\begin{array}{l}
\displaystyle {\bf v}^{(\infty)}\Big(u_{xt}-\alpha
u-\frac{1}{3}\beta
(u^3)_{xx}\Big)=0,\\[2mm]
\displaystyle u_{xt}=\alpha u+\frac{1}{3}\beta (u^3)_{xx},\\[2mm]
\displaystyle u_{x^2t}={\rm D}_x\big(\alpha u+\frac{1}{3}\beta (u^3)_{xx}\big),\\[2mm]
\displaystyle u_{xt^2}={\rm D}_t\big(\alpha u+\frac{1}{3}\beta (u^3)_{xx}\big),\\
\qquad \vdots\\ \displaystyle u_{xt^3}={\rm D}^2_t\big(\alpha
u+\frac{1}{3}\beta (u^3)_{xx}\big),
\end{array}
\end{eqnarray}
which leads to a polynomial of $u_{x^5}$ and $u_{t^5}$, with
functional coefficients of
\begin{eqnarray}
Q(x,t,u,u_x,u_t,u_{xx},u_{tt},,u_{x^3},u_{t^3},u_{x^4},u_{t^4})
\end{eqnarray}
and its derivatives. All of its coefficients must be zero. This
leads to a system of $5$ linear determining PDEs:
\begin{eqnarray}\label{eq:23}
&&\hspace{-10mm}\beta^4u_{xx}^8Q_{u_{t^4}^2}+Q_{u_{x^4},u_{t^4}}=0,\nonumber\\
&&\hspace{14mm}\vdots\\
&&\hspace{-10mm}u_xu_{tt}Q_{u,u_t}+\cdots+u_xu_{t^4}Q_{u,u_{t^3}}=0.\nonumber
\end{eqnarray}
Therefore, the most general third-order characteristic function
$Q$ is
\begin{eqnarray}\label{eq:24}
&&\hspace{-10mm}Q=\left(c_1\,t+c_2\right) u_t+3c_1u-c_1\,xu_x+c_3\,u_{t^3}-c_3\,\beta^3u_{xx}^6u_{x^3}\\
&&\hspace{-3mm}-\frac{3}{2}c_3\,\alpha\beta^2u_x
u_{xx}^4-\left(c_3\,\beta\alpha^2u_x^2-c_5\right)
u_x+{\frac{c_4\,u_{x^3}}{\sqrt {2\beta
u_{x^3}^2+\alpha}}},\nonumber
\end{eqnarray}
where $c_1,\cdots,c_5$ are arbitrary constants. There is not any
non-trivial second or forth-order characteristics. Thus, we prove
that
\paragraph{Theorem}
{\it The most general third-order infinitesimal local symmetry
generator of SPE (\ref{SPE}) is a ${\Bbb R}-$linear combination of
following five vector fields ${\bf v}_1$, ${\bf v}_2$, ${\bf v}_3$
of (\ref{eq:9}) and
\begin{eqnarray}\label{eq:25}
&&\hspace{-10mm}{\bf v}_4={\frac{\,u_{x^3}}{\sqrt {2\beta
u_{x^3}^2+\alpha}}}\,\partial_u,\nonumber\\[-4mm]
\\&&\hspace{-10mm}{\bf
v}_5=\Big(u_{x^3}-\beta^3u_{xx}^6u_{x^3}-\frac{3}{2}\alpha\beta^2u_xu_{xx}^4-\alpha^2\beta
u_x^3\Big)\,\partial_u.\nonumber
\end{eqnarray}
There is not any non-trivial second or forth-order infinitesimal
local symmetry generators.}

\end{document}